\documentclass[12pt]{article}
\usepackage{amscd, amssymb}
\usepackage{enumerate}
\usepackage[utf8]{inputenc}
\usepackage{hyperref}

\newcommand{\comment}[1]{}
\newcommand{\sr}{\rightarrow}

\newcommand{\nn}{{\bf N\rm}}

{\vskip
3mm}

\newcommand{\coq}[1]{``{\tt #1}"}
\begin{document}

\parskip = 2mm
\begin{center}
{\bf\Large An experimental library of formalized mathematics based on the univalent foundations.}

\vspace{3mm}

{\large\bf Vladimir Voevodsky}\footnote{School of Mathematics, Institute for Advanced Study,
Princeton NJ, USA. e-mail: vladimir@ias.edu}$^,$\footnote{Work on this paper was supported by NSF grant 1100938.}
\vspace {3mm}

\end{center}

\subsection{Introduction}
This is a short overview of an experimental library of mathematics formalized in the Coq proof assistant using the univalent interpretation of the underlying type theory of Coq.  I started to work on this library in February 2010 in order to gain experience with formalization of  mathematics in a constructive type theory based on the intuition gained from the univalent models (see \cite{KLV1}). 

Univalent models interpret types not as sets but as homotopy types. Their use in formalization of general mathematics (as opposed to just homotopy theory) is based on the following consideration. First note that we can stratify mathematical constructions by their ``level". There is element-level mathematics - the study of element-level objects such as numbers, polynomials or various series. Then one has set level mathematics - the study of sets with structures such as groups, rings etc. which are invariant under isomorphisms.  The next level is traditionally called category-level, but this is misleading. A collection of set-level objects naturally forms a groupoid since only isomorphisms are intrinsic to the objects one considers, while more general morphisms  can often be defined in a variety of ways. Thus the next level after the set-level is the groupoid-level - the study of properties of groupoids with structures which are invariant under the equivalences of groupoids. From this perspective a category is an example of a groupoid with structure which is rather similar to a partial ordering on a set.  

Extending this stratification we may further consider 2-groupoids with structures, n-groupoids with structures and $\infty$-groupoids with structures. Thus a proper language for formalization of mathematics should allow one to directly build and study groupoids of various levels and structures on them. 

A major advantage of this point of view is that unlike $\infty$-categories, which can be defined in many substantially different ways the world of $\infty$-groupoids is determined by Grothendieck correspondence (see \cite{Esquisse}) , which asserts that $\infty$-groupoids are ``the same" as homotopy types. Combining this correspondence with the previous considerations we come to the view that not only homotopy theory but {\em the whole of mathematics}  is the study of structures on homotopy types.  

The univalent models of constructive type theories enable one to use such type theories to reason directly about homotopy types with structures. This is the main idea of the Univalent Foundations of Mathematics - to use constructive type theory together with the intuition coming from its univalent homotopy-theoretic semantics to write and to prove theorems about mathematical objects of all ``levels" formally. 

Univalent Foundations can be seen as a realization of the vision of Michael Makkai whose paper \cite{Makkai} was very important for me in my search for a formal language for contemporary mathematics.

At the moment there are two actively supported proof assistants based on constructive type theories - Coq and Agda. Both proof assistants continue to be developed by teams which consist mainly of computer scientists who are actively experimenting with new features which are introduced into the systems without a formal verification of their consistency. Of these two systems Coq has been in development longer and is more conservative.  To further minimize the possibility of accidentally using a feature which may later be found to be inconsistent  the library described here was written using a restricted subset of the type theory underlying Coq. For another approach in Coq we suggest the reader look at the HoTT project library at \url{https://github.com/HoTT/HoTT}. The version of the library which this text refers to was checked to compile with a patched version of {\tt Coq 8.4pl3}. For instructions on how to get this version of Coq and how to patch it see the file {\tt Coq\_patch/README}. 

The type theory of Coq is, roughly speaking, a combination of three components. The first component is a version of the Thierry Coquand's Calculus of Constructions (CC) (see \cite{Coquand}). This is a type system with two universes,  Prop and Type, dependent products and abstraction/application constructions satisfying $\beta$-reduction. The second component is a universe management system which replaces two universes of CC with an infinite hierarchy of universes which is due to Z. Luo (see \cite{Luo}). The third component is a machinery for creating strictly positive ``inductive types" described in  \cite{Mohring1}. 

In our library we use a small subset of a modified version of the Coq type system. The modifications are introduced through a patch contained in the subdirectory Coq\_patch. Some information on the content of this patch and on its history can be found in the README file of that subdirectory. 

The main modification turns off the universe consistency verification system of Coq. This of course makes the type system inconsistent (any type, including the empty type, can be shown to have an object). The proper solution is instead to use universe polymorphism together with either resizing rules (see \cite{VUFresizingslides}) or higher inductive types (see \cite{hottbook}). However these modifications are highly non-trivial and for the experimental purposes of the current library it seemed reasonable to rely on careful tracing of universe levels ``by hand". This issue becomes important only starting with the file {\tt hProp.v}. The first major file of the library, {\tt uu0.v}, can be compiled without the patch. 

The main restriction which we impose on the constructions of the library concerns the use of the inductive types machinery of Coq. In a rather ingenious way this machinery is normally used in Coq both to define standard ingredients of constructive type theories such as identity types, dependent sums, the one point type, disjoint unions, the empty type,  booleans and natural numbers and also to define a multitude of other constructs such as, for example, inequalities between natural numbers. In the current library we use this machinery only to introduce the standard constructions listed above. No further use of inductive types is made except in one place in the file {\tt hnat.v} where we show that our approach to comparisons between natural numbers is equivalent to the approach  taken in the standard library of Coq. 

Another restriction is that we do not use the universe Prop. Associated with this universe there is a ``singleton elimination" rule which is inconsistent with the univalent model. To avoid accidental use of this rule by tactics the patch file modifies the way the universe level of inductive constructions (most notably of identity types) is computed. During the compilation of the first file of the library, {\tt uuu.v}, the compiler should display \coq{paths 0 0:UUU}. Without proper application of  the patch the compiler would display \coq{paths 0 0:Prop}. 

The distribution of Coq includes an extensive ``standard library". Our library uses only the first and most basic subdivision of the Coq's standard library, namely Coq.Init. In fact some of the files of the standard library may take very long to compile with the "-no-sharing" option which is introduced by the patch and which we use to overcome a bug in Coq's normalization algorithm. See the file {\tt Coq\_patch/README} for instructions of how to compile Coq without compiling most of the standard library.

\subsection{File {\tt uuu.v}}

The first several lines of {\tt uuu.v} introduce new notations for some of the constructions that are defined in Coq's standard library. The part of the library where these constructions are introduced is  located in \coq{coqlocation/theories/Init/} where \coq{coqlocation} is the directory where the Coq distribution is. Files of this part of the library are automatically loaded by Coq while to load other parts of the standard library (located in other subdirectories of \coq{coqlocation/theories/}) requires an explicit instruction.  

In the first new definition in \coq{uuu.v} we introduce the version of {\em dependent sum} used in our library. It is called ``{\tt total2}" due on the one hand to its semantic meaning as the total space of a fibration and on the other to its function as a generic record of length 2\footnote{In the first version of the library there was also \coq{total3} corresponding to the generic record of length 3.}. Several important features of Coq formalization can be illustrated with this definition and the following definitions of ``{\tt pr1}" and ``{\tt pr2}".

The first parameter of the construction, the type ``{\tt T}", is shown in the definition in braces. This means that this is an {\em implicit} parameter, i.e., when ``{\tt total2}" is used one writes \coq{total2 P} instead of \coq{total2 T P}. Types of expressions are computable in Coq from expressions themselves and since the type of ``{\tt P}" must be ``{\tt T->Type}" the system can infer ``{\tt T}" from ``{\tt P}". 

The second parameter ``{\tt P}" is of the type ``{\tt T->Type}". Here ``{\tt Type}" is a generic notation which Coq uses for universes. The universe management in Coq is rather baroque and well hidden from user control so for simplicity one may think that ``{\tt Type}" is synonymous with the name of some fixed universe ``{\tt UU}". A function ``{\tt T->UU}" is intuitively a way to assign to any object of ``{\tt T}" an object of ``{\tt UU}", i.e., a type which is contained in ``{\tt UU}". In other words, it is a family of types in ``{\tt UU}" parametrized by ``{\tt T}". The semantics of this is as follows. If we use the univalent model with values in the category of simplicial sets then ``{\tt T}" is mapped to a (Kan) simplicial set and ``{\tt UU}" is mapped to the base of the universal Kan fibration which classifies Kan fibrations whose fibers belong to ``{\tt UU}".  Thus ``{\tt P}" corresponds to a Kan fibration over ``{\tt T}" and ``{\tt total2 P}" is the total space of this fibration. 

In the informal semantics with values in $\infty$-groupoids ``{\tt T}" is mapped to an $\infty$-groupoid while ``{\tt UU}" is mapped to the $\infty$-groupoid of $\infty$-groupoids in ``{\tt UU}" and their equivalences. The function ``{\tt P}" then  can be viewed as a functor and ``{\tt total2 P}" is the $\infty$-groupoid of pairs $(x, y)$ where $x$ is an object of ``{\tt T}" and $y$ an object of $P(x)$. 

The next definition is that of ``{\tt pr1}". It takes three parameters and returns an object of ``{\tt T}" where ``{\tt T}" is the first parameter. In general, when one has a definition of some ``{\tt C}" in Coq with parameters ``{\tt x1 x2... xn}" one can write not only ``{\tt C a1... an}" but also {\em partially applied} versions such as ``{\tt C a1... a(n-1)}" or just ``{\tt C}". The type of such a partially applied definition will be a function type or more generally a {\em dependent product} type. 

In the case of ``{\tt pr1}" the first two parameters are implicit and supposed to be inferred from the third. If one wants to use a partially applied version of ``{\tt pr1}" one has to provide the first two parameters explicitly. To tell Coq that a definition will be used as if all its parameters were explicit one uses prefix ``{\tt @}" and writes, for example, ``{\tt @pr1 T P}". The type of this expression is the function type ``{\tt (total2 T P) -> T}" and its semantical meaning is the projection from the total space of a fibration to its base.

An extremely important feature of dependent type theories which is unavailable in the theories without dependent types and which at the first may seem confusing is that we also have ``{\tt @pr2 T P}". Obviously not any fibration is trivial so we do not normally have a projection from the total space to a fiber {\em as a function}. However we always have it as a {\em dependent function}. By writing something like

{\tt Variable T:Type.

Variable P:T -> Type.

Check ( @pr2 T P ).}

one will see that the type of  ``{\tt @pr2 T P}" is ``{\tt forall tp:total2 P, P (pr1 tp)}".  The semantic meaning of the later expression is as follows. ``{\tt forall}" is the name of the dependent product construction in Coq. Its general format is ``{\tt forall x:T1, T2}" where ``{\tt T1}" is a type expression and ``{\tt T2}" is a type  expression which may have a parameter ``{\tt x}" of type ``{\tt T1}".  Such an expression with a parameter semantically is the same as a function ``{\tt T1 -> UU}", i.e., the ``{\tt forall}" construction has essentially the same parameters as ``{\tt total2}" - a type and a family of types parametrized by this type. As was explained above such a pair corresponds in the univalent model in simplicial sets to a Kan fibration. The type   ``{\tt forall x:T1, T2}" is the (Kan) simplicial set of {\em sections} of this fibration. 

If ``{\tt T2}" does not actually depend on ``{\tt x}" then one abbreviates the expression  ``{\tt forall x:T1, T2}" to ``{\tt T1 -> T2}". Semantically it corresponds to the case of a constant fibration whose sections are just functions from the base ``{\tt T1}" to the fiber ``{\tt T2}". 

Returning to the case of ``{\tt @pr2 T P}" we see that semantically it is a section of the fibration over ``{\tt @total2 T P}" whose fiber over ``{\tt tp}" is the fiber of ``{\tt P}" over ``{\tt pr1 tp}". In mathematical notation, if our fibration is $p:E\sr B$ then ``{\tt @pr2 T P}" is the diagonal section of $E\times_B E$ over $E$. 

\subsection{File {\tt uu0.v}}

This file contains the results of the library which are applicable to {\em all} types. 

The first three lines of the file are also repeated with some obvious changes in all the rest of the files of the library. These are commands to the Coq program. 

The first one tells Coq not to do a certain type of steps automatically at the start of every proof but to leave the choice of whether or not to do these steps to the user.

The second and the third lines address the mechanism which loads other library files. They are discussed in more detail in the Appendix. 

Let me now use some of the first proofs given in \coq{uu0.v} to illustrate how the proof system of Coq works.  Note first that a line such as 

``{\tt Definition name1:expr1.}"

tells Coq that a constant called ``{\tt name1}" of type ``{\tt expr1}" will be provided by the user. In the case of ``{\tt Definition}" there are two ways to provide the value of this constant. One can write

``{\tt Definition name1:expr1:= expr2.}"

in which case ``{\tt expr2}" should be an expression which has type ``{\tt expr1}" which will be the value of the constant ``{\tt name1}". Alternatively, one can write ``{\tt Proof.}" after ``{\tt Definition name1:expr1.}" and then use various commands of Coq proof mode to construct the value of the constant. When Coq says ``{\tt Proof completed}" in the ``{\tt response}" window one writes either ``{\tt Qed.}" or ``{\tt Defined.}". The difference between the two is that when ``{\tt Qed.}" is used the actual structure of the constructed expression becomes hidden (opaque) while when ``{\tt Defined.}" is used the structure remains accessible.  

The keywords ``{\tt Theorem}", ``{\tt Lemma}" and ``{\tt Proposition}" are strictly equivalent and are equivalent to ``{\tt Definition}" except that one must use the proof mode to provide the value of the corresponding constant, i.e., one can not simply provide the value after ``{\tt :=}". 

More generally Coq can be told that a constant with the name ``{\tt name1}" is going to be introduced by a line of the form

``{\tt Definition name1 ( x1:texpr1 )... ( xn:texprn ):expr.}"

which is essentially equivalent to 

``{\tt Definition name1:forall x1:texpr1,..., forall xn:texprn, expr}"

with the only difference being that the first form allows one to say that some of the parameters will be implicit by using curly brackets.  

The first proof of the library is that of ``{\tt Definition fromempty}". The sentence which starts with the word ``{\tt Definition}" tells Coq that a constant with the name ``{\tt fromempty}" of type ``{\tt forall X:UU, empty -> X}" will be provided and that the type parameter ``{\tt X}" is implicit. The value for this constant is constructed inside the proof mode through the use of two tactics ``{\tt intros}" and ``{\tt destruct}". We will not discuss here how the tactics language of Coq is working referring the reader instead to  Coq Reference Manual.

Detailed information about the mathematical content of the file {\tt uu0.v} can be obtained from the comments in this file. We will only discuss here a few fundamental constructions the meaning of which might not be immediately obvious.

The first such construction is ``{\tt iscontr T}" where ``{\tt T}" is a type. It introduces the concept from which almost everything else is build - the concept of a contractible type. By definition, a proof of contractibility of a type ``{\tt T}" is an object of the type ``{\tt iscontr T}". There are two ways to argue that this is a ``{\tt correct}" way to define contractibility. The first one is to point out that the more complex homotopy-theoretic notions defined with the use of this notion of contractibility are proved further on in this file to satisfy a large number of expected properties.   

Another is to analyze the univalent semantics of this construction.  Consider for example a univalent model with values in Kan simplicial sets. Then ``{\tt T}" is a simplicial set. A point of ``{\tt iscontr T}" is a pair ``{\tt (cntr, s)}" where ``{\tt cntr}" is a point of ``{\tt T}" and ``{\tt s}" is an object of ``{\tt forall t:T, paths  cntr t}". The family of types ``{\tt t $\mapsto$ paths cntr t}" is the paths bundle corresponding to the point ``{\tt cntr}" and as explained above ``{\tt s}" is a section of this bundle. But the total space of the paths bundle is contractible and if it has a section then ``{\tt T}" is a retract of a contractible simplicial set and therefore it is contractible. In the opposite direction if ``{\tt T}" is contractible then it is in particular non-empty and we can choose a point ``{\tt cntr}" in ``{\tt T}". Any fibration over a contractible s.s. is trivial and if it has a non-empty fiber it has a section. The fiber of the paths fibration defined by ``{\tt cntr}" over ``{\tt cntr}" is non-empty and therefore it has a section ``{\tt s}" which gives us a point   in ``{\tt iscontr T}". 

The next fundamental definition is the property ``{\tt isweq}" of a function ``{\tt f}" to be a (weak) equivalence which is defined as the condition that all (homotopy) fibers of ``{\tt f}" are contractible. Along with ``{\tt isweq f}" we introduce ``{\tt weq X Y}" - the type of (weak) equivalences from ``{\tt X}" to ``{\tt Y}", i.e., of pairs ``{\tt (f, is)}" where ``{\tt f:X -> Y}" and ``{\tt is:isweq f}".  

Theorem ``{\tt gradth}" shows that for a homotopy equivalence, i.e., a quadruple ``{\tt f:X -> Y}", ``{\tt g:Y -> X}", ``{\tt egf}", ``{\tt efg}" where ``{\tt egf}" is a homotopy from \coq{funcomp f g} to the identity of ``{\tt X}" and  ``{\tt efg}" is a homotopy from \coq{funcomp g f} to the identity of ``{\tt Y}", the function ``{\tt f}" is a (weak) equivalence. The difference between the notions of a homotopy equivalence and a weak equivalence is somewhat subtle but important. Let ``{\tt X}" and ``{\tt Y}" be types and ``{\tt homeq X Y}" the type of quadruples ``{\tt (f, (g, (egf, efg)))}". Theorem ``{\tt gradth}" (or rather Definition ``{\tt weqgradth}") defines a function ``{\tt (homeq X Y) -> (weq X Y)}".  Using Definitions ``{\tt homotweqinvweq}" and ``{\tt homotinvweqweq}" one gets a function ``{\tt (weq X Y) -> (homes X Y)}" and it is not difficult to show that these functions make ``{\tt weq X Y}" into a retract of ``{\tt homeq X Y}". In general however this retraction is not an equivalence.  The reason why \coq{weq} is ``better" than \coq{homeq} is related to the difference between  {\em properties} and {\em structures} which is explained below.

Corollary \coq{iscontrweqf} and Definition \coq{wequnittocontr} show that a type is contractible iff it is weakly equivalent to \coq{unit} and in particular that up to weak equivalence there is only one contractible type. Corollary \coq{isweqmaponpaths} shows that a weak equivalence defines a weak equivalence on \coq{paths} types. Theorems \coq{twooutof3a}, \coq{twooutof3b} and \coq{twooutof3c} establish the 2-out-of-3 property of weak equivalences - if two out of three functions \coq{f}, \coq{g}, \coq{funcomp f g} are weak equivalences then so is the third. All these results are proved using \coq{gradth}.

Then  there follows a series of simple results which assert that various natural functions such as the ones defining associativity and commutativity of direct products or distributivity of direct products and binary coproducts are weak equivalences.

The next tool-box which we introduce contains the type-theoretic versions of the results and definitions related to homotopy fiber sequences.  Our approach to fiber sequences differs somewhat from the usual approaches. A fiber sequence structure \coq{fibseqstr} on a triple \coq{f:X->Y, g:Y->Z, z:Z} is defined as a homotopy from \coq{funcomp f g}  to the constant function \coq{fun x : X => z} such that the associated function \coq{ezmap} from \coq{X} to the homotopy fiber \coq{hfiber f z} is a weak equivalence.  

For any fiber sequence structure \coq{fs} on \coq{(f,g,z)} and any object \coq{y:Y} we construct a function \coq{d1: paths (g y) z -> X} and the {\em derived} fiber sequence structure \coq{fibseq1} on the triple \coq{(d1, f, y)}. This construction can be iterated leading to a type theoretic construction of long homotopy exact sequences of fibrations.

We then investigate three standard situations where fibre sequences arise.

For any family of types \coq{P:Z->UU} over a type \coq{Z} and an object \coq{z:Z} we construct in \coq{fibseqpr1} a fiber sequence structure on the triple \coq{(iz,pr1,z)} where \coq{iz} is the inclusion of the fiber \coq{P z} to \coq{total2 P} and \coq{pr1} the projection \coq{total2 P -> Z}. Applying to it the construction of the derived fiber sequence we get a family of weak equivalence \coq{ezweq1pr1} which connect the homotopy fibers of \coq{iz} with paths types on \coq{Z}.

For a function \coq{g:Y->Z} and an object \coq{z:Z} we define in \coq{fibseqg} the obvious structure of a fiber sequence on the triple \coq{(hfiberpr1, g, z)} where \coq{hfiberpr1 : hfiber g z -> Y} is the standard function and give explicit descriptions of its first, second and third derived sequences.

Finally we construct a fiber sequence \coq{fibseqhf} of homotopy fibers of a composable pair of functions \coq{f:X->Y}, \coq{g:Y->Z} for \coq{z:Z} and  \coq{ye: hfiber g z} with the underlying sequence of morphisms of the form \coq{hfiber f (pr1 ye) -> hfiber (comp g f) z -> hfiber g z}.

The next fundamental notion which we introduce is the notion of h-levels. The definition ``{\tt isofhlevel n}" uses the type \coq{nat} of natural numbers which is introduced in Coq.Init.Datatypes as the inductive type with two constructors \coq{O} of type \coq{nat} (corresponding to $0$) and \coq{S} of type \coq{nat -> nat} (corresponding to the successor function $n\mapsto n+1$).  Semantically we have that $T$ is of h-level $0$ iff it is contractible and of h-level $1+n$ iff for any $x,y$ in $T$ the paths space $paths\,x\,y$ is of h-level $n$.

A function \coq{f:X->Y} is said to be of h-level \coq{n} if all its (homotopy) fibers are of h-level \coq{n}. In particular, a function is of h-level \coq{O} iff it is a weak equivalence.  

Types of h-level 1 are called {\em propositions} and we write \coq{isaprop} instead of \coq{isofhlevel 1}. A homotopy type $T$ is of h-level $1$ iff for any $x,y\in T$ the paths space between $x$ and $y$ is contractible.  In the world of classical homotopy types there are  only two homotopy types with this property - the empty type and the contractible type. If \coq{T} is of h-level 1 and it is inhabited, i.e., there is an object \coq{t:T} then, as \coq{iscontraprop1} shows \coq{T} is contractible.  However there are many non-equivalent types of h-level 1 which have no objects. This discrepancy between the model side and the syntactic side is the univalent form of the first Goedel's incompleteness theorem. 

It is of a  fundamental importance for the univalent approach to distinguish types which are propositions from more general types. In particular, if one wants to formalize univalently non-constructive proofs then one should add the axiom of excluded middle to the environment. Adding it in the form \coq{forall T:Type, coprod T (T-> empty)} would be incompatible with the univalent models (and with the univalence axiom). This however does not mean that univalent semantics is incompatible with classical logic - the correct univalent formulation of the theorem of excluded middle is  \coq{forall T:hProp, coprod T (T-> empty)} where \coq{hProp := total2 (fun T:Type => isaprop T)}. 

A function between classical homotopy types $f:X\sr Y$ is of h-level 1 iff it is homotopy equivalent to the inclusion of a union of connected components of $Y$ into $Y$. On the type theoretic side we define inclusions as functions of h-level 1 (\coq{isincl}). 

Inclusions correspond to {\em predicates} or {\em properties} - functions \coq{P:T->UU} such that \coq{forall t:T, isaprop (P t)}. Given such \coq{P} we can form the type  \coq{total2 P} whose objects are pairs \coq{(x,p)} where \coq{x:T} and \coq{p:P x}. By \coq{isweqezmappr1} the homotopy fiber of the projection \coq{pr1:total2 P -> T} over \coq{x:T} is weakly equivalent to \coq{P x}. By \coq{isofhlevelweqf} the h-levels are invariant under weak equivalences. We conclude that the projection \coq{total2 P -> T} is a (homotopy) inclusion iff for all \coq{x:T} the type \coq{P x} is of h-level $1$. If the h-level of \coq{P x} is greater than $1$ for some \coq{x:T} then \coq{P} defines a {\em structure} on objects of \coq{T}. 

One of the important naming conventions in our library is that any name which starts with ``is" such as \coq{isontr} or \coq{isweq} corresponds to a {\em property}.  For further discussion of propositions and properties in the univalent approach see Section \ref{hProp}.

Types of h-level 2 are called sets (or, sometimes, h-sets) and we write \coq{isaset} instead of \coq{isofhlevel 2}. A classical homotopy type $T$ is of h-level 2 iff the path space between any two points is either empty or contractible -  one can easily see that this is equivalent to the condition that $T$ is a disjoint union of contractible components, i.e., that it is homotopy equivalent to a set. On the type-theoretic side, due to the constructive nature of the theory, sets need not be disjoint unions of points. More precisely it is not necessarily true that for a set \coq{T} and an object \coq{t:T} there is an equivalence between \coq{T} and \coq{coprod  (compl T t) unit} where \coq{compl T t} is the complement to \coq{t} in \coq{T}.  Types which satisfy the later property for all objects are called {\em types with decidable equality} (see \coq{isdeceq}). We show that any type with decidable equality is an h-set in \coq{isasetifdeceq} and use it to prove that Booleans (\coq{isasetbool}) and natural numbers (\coq{isasetnat}) are h-sets but not all h-sets can be proved to have decidable equality.  A simple example of an h-set which does not have decidable equality is the type of functions \coq{nat -> Bool}. Such types as Dedekind reals or p-adic numbers are also h-sets with undecidable equality.   

Most of mathematics as we know it deals with structures of h-level 2 on types of h-level 2. For example, a group is a pair \coq{(T,S)} where \coq{T} is an h-set and \coq{S} is an object of the h-set of group structures on \coq{T}. For further discussion of h-sets see Section \ref{hSet}. 

For higher $n$ the notion of h-level coincides with the well known notion of n-types up to a shift of index by $2$, i.e., a type $T$ is of h-level $n+2$ iff for any $x$ in $T$ and $i>n$ one has $\pi_i(T,x)=0$. The best known area of mathematics whose univalent formalization requires types of h-level $3$ is {\em category theory}. For a univalent approach to category theory see \cite{RezkCompletion}.  

The file {\tt uu0.v} contains three axioms - \coq{funextempty}, \coq{etacorrection} and \coq{funextfunax} and the third one implies both the first and the second. Axioms are generally undesirable in constructive type theory even if, as is the case for these three axioms, they are semantically justified. The reason is that they tend to break a very important property of constructive type theories which is called {\em canonicity}. In its simplest form canonicity asserts that any object \coq{o} of type \coq{nat} (natural numbers) in the empty context which is in the {\em normal form} is of the form \coq{S ... S O}, i.e., is a {\em numeral}. We will come back to this property in the discussion of the files {\tt finitesets.v}, {\tt hz.v} and {\tt hq.v}. 

For example, \coq{funextfunax} can be used to define an object of type \coq{nat} which is in the normal form but which is not a numeral as follows. Consider the transport along a path \coq{transportf}. Let \coq{T} be any type constant defined in the empty context (e.g. \coq{unit} or even \coq{empty}). Let \coq{f:=fun t:T => t} be the identity function on \coq{T}. Let  \coq{e: paths f f} be the path obtained by applying \coq{funextfunax} to the homotopy \coq{fun t:T => idpath t t}. Set \coq{x := transportf (fun g:T -> T => nat) e O}. By doing this construction in Coq and typing \coq{Eval Compute in x} - the command which displays the normal form of expression \coq{x} - one immediately sees that \coq{x} does not normalize to a numeral. For a further discussion of this phenomenon and its relation to the problem of constructive interpretation of the univalence axiom see Section \ref{finitesets}.

Note that while an object of type \coq{nat} defined with the use of axioms {\em may} happen not to normalize to a numeral it is not necessarily so. In particular many of the test computations in files {\tt finitesets.v}, {\tt hz.v} and {\tt hq.v}  use theorems and definitions which include \coq{funextfunax}. 

Axiom \coq{funextfunax} is known as the {\em functional extensionality} axiom. In its original form it is not even an "axiom", i.e.,  the type of \coq{funextfunax} can not be proved to be a proposition. More precisely, one can show that the homotopy type corresponding to the type of \coq{funextfunax} under a univalent model has more than one connected component. To deal with this issue we use everywhere not \coq{funextfunax} itself but its corollary \coq{funcontr} which can be shown to be a proposition.  

In the following parts of the library we use \coq{funcontr} to show that the dependent product construction interacts in the expected way with weak equivalences (see \coq{isweqmaponsec} and \coq{isweqmaponsec1}) and with h-levels (see \coq{impred}).   We also prove a number of results which justify our use of ``is" prefix in the names of constructions such as \coq{iscontr}, \coq{isweq} and \coq{isofhlevel} by showing that the types of corresponding constants are indeed of h-level 1.

\comment{
\begin{description}
\item  [isweqmaponpaths] Weak equivalence on paths - if ``{\tt f:X -> Y}" is a weak equivalence then for any ``{\tt x1 x2:X}" the function ``{\tt maponpaths f:paths x1 x2 -> paths (f x1) (f x2)}" is a weak equivalence.
\item [twooutof3a, twooutof3b, twooutof3c] The 2-out-of-3 property - if two out of three functions ``{\tt f, g, g$\circ$ f}" are weak equivalences then so is the third.
\item [isweqcoprodf] Coproduct of two weak equivalences is a weak equivalence. 
\end{description}
}

\subsection{File {\tt hProp.v}}
\label{hProp}
This is the only (so far) file in the folder \coq{hlevel1}. It contains basic results related to types of h-level 1, i.e., to propositions. First we introduce the type \coq{hProp} which relates to types of h-level 1 in the same way as the universe \coq{UU} relates to all types. In fact we should consider the universe \coq{UU} as a parameter of \coq{hProp} writing \coq{hProp UU} for the type of propositions in a universe \coq{UU}.  Unfortunately the universe management system does not allow universe parameters and we are forced to consider \coq{hProp} with respect to a fixed universe \coq{UU}.

In the univalent model a type is a proposition iff it is empty or contractible. Therefore the model of \coq{hProp UU} is the simplicial subset of the model of \coq{UU} which consists of two connected components - the component of the empty type which is a 1-point simplicial set and the component of contractible types which is a large (relative to \coq{UU}) contractible simplicial set.   

This creates problems with the next construction in \coq{hProp} which we call \coq{ishinh\_UU} and which is also known as the bracket type or squash type construction. The idea is that for a type \coq{T} there should be a proposition \coq{ishinh\_UU T} which is true iff \coq{T} is inhabited. This is equivalent to saying  that  \coq{ishinh\_UU T} should be defined together with a function \coq{hinhpr T : T -> ishinh\_UU T} which is universal among functions from \coq{T} to propositions. Using   the fact that h-levels are stable under the formation of dependent products (\coq{impred} from {\tt uu0.v}) we show in \coq{isapropisinh} that \coq{ishinh\_UU T} is indeed a proposition and in in \coq{hinhuniv} that the function \coq{hinhpr} it is universal.

However there is an element of cheating here. In fact this part of {\tt hProp.v} would not go through in un-patched Coq. The only reason it works in Coq is that we use the patched version which does not check universe consistency. 

The problem is that \coq{ishinh\_UU T} is a proposition in a bigger universe than \coq{UU} which is universal with respect to functions from \coq{T} to propositions in \coq{UU}.

How can this problem be fixed without introducing potential inconsistency? There are currently three ideas. The first two have to do with {\em resizing rules} and the third with {\em higher inductive types}. All three are associated with interesting unsolved problems. Note that the issue is particular to the constructive setting. If we did not care about computation and added the excluded middle axiom then we could use a double negation version \coq{isinhdneg} of \coq{ishinh} which does not lead to any issues with universe levels.  

In the following part of the library we define an interpretation of intuitionistic logic on \coq{hProp}. The construction of \coq{ishinh\_UU} is a necessary prerequisite for the construction of the disjunction - the disjoint union of two propositions considered as types is not in general a proposition and one has to apply \coq{ishinh\_UU} to obtain disjunction as an operation on propositions.  

In the last part of the file we introduce the univalence axiom for \coq{hProp}  and consider some of its corollaries.

\subsection{File {\tt hSet.v}}
\label{hSet}

This file contains basic results related to sets, i.e., types of h-level 2.  The first brief section discusses types which satisfy axiom of choice, i.e., which are ``projective objects''. It is later used in {\tt stnfsets.v} and {\tt fintesets.v} to show that the axiom of choice holds for families over finite sets. 

Then we introduce a series of definitions and results about relations on types. Many of these results are later used to prove standard properties of comparisons on natural numbers and later on integers and rational numbers. 

The most important part of this file deals with set-quotients of types. The theory of quotients is well known to be one of the difficult points of the usual constructive type theory.  The univalent model provides an explanation for this fact - since types are homotopy types rather than sets the quotients need to be understood as homotopy quotients which are often very complicated. 

The set-quotients are, from homotopy-theoretical point of view, quotients with respect to ``homotopy-invariant equivalence relations". The finest such relation is given by the condition ``$a$ is path-connected to $b$'' with the corresponding quotient being $\pi_0$. Quotients with respect to stronger equivalence relations on $T$ are quotients of $\pi_0(T)$. The quotient with respect to the strongest relation, i.e., the one where any two points are equal is equivalent to \coq{ishinh\_UU T}. 

While in classical setting such quotients create no problems in the constructive setting things are more complicated. One problem is the increase in the universe level when one passes to a quotient. It is similar to the  problem which we discussed in the context of \coq{ishinh\_UU}. 

Another problem can be seen in the way in which taking quotients interacts with taking sub-objects. Let \coq{X} be a type, \coq{R} an equivalence relation on \coq{X} and \coq{P:setquot R -> hProp} a predicate on the quotient of \coq{X} with respect to \coq{R}. The composition \coq{Q} of \coq{P} with the projection \coq{setquotpr R:X -> setquot R} is a predicate  on \coq{X}.  Let \coq{U:= carrier P} and \coq{X':= carrier Q} be they sub-objects of \coq{setquot R} and \coq{X} respectively corresponding to \coq{P} and \coq{Q}. The restriction \coq{R'} of \coq{R} to \coq{X'} is an equivalence relation and we may consider two types \coq{setquot R'} and \coq{U}. As is proved in \coq{weqsubquot} these two types are equivalent. However the equivalence \coq{U -> setquot R'}, the obvious function \coq{X' -> U} and the projection \coq{setquotpr R': X' -> setquot R'} do not commute {\em computationally}. 

This leads for example to the use of somewhat unnatural constructions to define the inverse on non-zero elements of fields of fractions since the straightforward definition ``does not compute". 

A possible way to deal with this issue by extending the type theory of Coq with a new component called tfc-terms (from the trivial fibration/cofibration axiom of model categories) is briefly discussed in the comments after \coq{weqsubquot}. 

At the end of the file {\tt hSet.v} we describe another approach to set-quotients. Originally this part was written because I thought that the computational behavior of this alternative construction will be better. However it turned out to have very similar (and probably equivalent) problems as the first one.

\subsection{Files {algebra1*.v}}

These files introduce the standard notions of abstract algebra including the interaction between algebraic operations and partial orderings. 

The file {\tt algebra1a.v} introduces basic definitions related to binary operations and pairs of binary operations on h-sets. A few definitions where the generalizations from h-sets to all types are straightforward are given for arbitrary types. 

The file {\tt algebra1b.v} is about monoids, abelian monoids, groups and abelian groups including the construction of monoids of fractions in the abelian case. 

The file {\tt algebra1c.v} is about rigs (semi-rings with a unit such that $0\cdot 1=1\cdot 0 = 0$), commutative rigs, rings and commutative rings. It includes the construction of the ring of differences from a rig and of localization of a commutative ring by a multiplicative system of elements.   We also prove the basic results about the behavior of partial orderings and equivalence relations with respect to these constructions.  

The file {\tt algebra1d.v} is the first one which contains material which is probably unusual for an average mathematician. It deals with the notions of an integral domain and of a field in constructive framework. Unlike the notions considered above the notions of an integral domain and of a  field acquire additional distinctions in constructive mathematics relative to the classical one. For example the condition ``every non-zero element is invertible" in the definition of a field has three non-equivalent constructive formulations - one can require that any element which is non-invertible is zero or that any element which is non-zero is invertible or that any element is either invertible or equals zero. 

In {\tt algebra1d.v} we consider the later definition (any element is either invertible or equals zero). It is the strongest (most restrictive) one and it immediately implies that the equality on a field is a decidable relation. This is clearly unsatisfactory for many purposes - for example real numbers or the ``field" of power series do not satisfy this condition.  To deal with this problem one needs to introduce the notion of apartness relations and study their interactions with algebraic structures. Some information on the subject as well as further formalizations in the style of this library can be found in \cite{PVW1}. 

In {\tt algebra1d.v} we restrict ourselves to the case of decidable equality and give in that case a constructive definition of a field of functions of a (decidable) integral domain. 

All constructions in the algebra files have non-trivial extensions from h-sets to arbitrary types. For example the notion of a monoid generalizes to as yet undefined notion of an H-type which should include all the higher coherence structures associated with associativities. The notion of a partially ordered set generalizes to the notion of $(\infty,1)$-category and the notion of a partially ordered monoid generalizes to the notion of a monodical $(\infty,1)$-category.  I do not know what is the classical name for the higher analogs of rigs and commutative rigs (going from rigs to rings is straightforward since the only axiom involved is the invertibility of addition which has a formulation common for types of all levels) and whether such objects been considered.  None of these have as yet been defined in terms of type theory.

\subsection{File {\tt hnat.v}} 

In this file we provide basic constructions and results related to the arithmetic operations and comparisons on natural numbers. The type \coq{nat} is introduced in Coq.Init. We use this definition for natural numbers and also standard definitions for the addition, subtraction (which is defined such that for $n<m$ one has $n-m=0$) and multiplication on \coq{nat}. 

Our approach to comparisons is different from the one used in Coq.Init. There the main comparison is \coq{le} which is introduced through an inductive definition based on the principle that \coq{le} is a family of types whose objects are either "reflexivity" comparisons in \coq{le n n} or successor comparisons obtained from constructor of the form \coq{le n m -> le n (S m)}. 

Since our library uses only those inductive constructions in Coq which are necessary for the definition of the standard ingredients of the Martin-L{\"o}f Type Theory we do not use \coq{le}. Instead we start with Boolean "greater" which we call \coq{natgtb} defined by induction on \coq{nat} as a function \coq{nat -> nat -> Bool}, define \coq{natgth n m} as \coq{paths (natgtb n m) true} and then define the three other comparisons \coq{natlth}, \coq{natleh} and \coq{natgeh} in terms of \coq{natgth}. This has the advantage that the same definitions of ``less", ``less or equal" and ``greater or equal" in terms of ``greater" work for integers and rationals and the proofs of the main properties of these comparisons from the main properties of ``greater" can be directly copied from the \coq{nat} case to the cases of \coq{hz} and \coq{hq}. 

After this choice of how to define the comparisons and prove their properties is made the rest is rather straightforward. 

At the end of the file we analyze the Coq.Init construction of \coq{le}-types showing that \coq{le n m} is always a proposition (i.e., has h-level 1). 

\subsection{File {\tt stnfsets.v}} 

This is the first of the two files where we introduce constructions related to finite sets. In this file we deal only with ``standard" finite sets which are defined such that \coq{stn n} is the type of natural numbers which are less than \coq{n}. 

Most of the file is occupied by constructions of various weak equivalences involving standard finite sets. For example we construct a weak equivalence between \coq{weq (stn n) (stn n)} and \coq{stn (factorial n)}. 

At the end of the file we use the notion of a standard finite set to formulate and prove results on bounded quantification and then to give a univalent proof of the accessibility theorem for natural numbers. 

\subsection{File {\tt finitesets.v}}
\label{finitesets}

We define the structure of having $n$ elements on a type \coq{T} as a weak equivalence from the standard set with $n$ elements to \coq{T}. A type \coq{T} is called a finite set if there exists (or, in terminology of \cite{hottbook}, if there merely exists) a pair \coq{tpair \_ n s} where \coq{n:nat} and \coq{s} is a structure of having $n$ elements on \coq{T}.  We then use the results of {\tt stnfsets.v} to show that various constructions on finite sets produce finite sets. 

An important property of our approach is that despite the fact that we use ``mere" existence in the definition of what it means to be finite there is a function \coq{fincard} which computes the cardinality of a finite set. 

Related to this function are several examples of computation which are included at the end of the file {\tt finitesets.v}. The property of Martin-L{\"o}f type theory which makes automatic computation possible is known as {\em canonicity}. In its simplest form the {\em canonicity theorem} asserts that 
%
%???? - explain normal form, normalization/computation somewhere?
%a
any object \coq{o} of type \coq{nat} defined in the empty context which is in the normal form is  a numeral, i.e., a sequence  \coq{S ... S O} (recall that \coq{O} is the notation for $0\in \nn$ and \coq{S} is the notation for the successor function $n\mapsto 1+n$).  

The possibility of automatic terminating computation is a corollary of this property combined with {\em strong normalization} - the assertion that any sequence of reductions starting with a given well-formed expression is finite\footnote{Strong normalization is a difficult theorem. In particular, using a variant of Goedel's argument, it can be shown that it can not be proved unconditionally. In practice, all known proofs of strong normalization for Martin-L{\"o}f type theory require one to assume that a substantial portion of ZFC is consistent.}. 

By definition, an expression is said to be in the normal form if there are no reduction steps starting with this expression. Therefore in a theory with strong normalization for any well-formed expression there is a finite sequence of reductions which results in an expression in the normal form. If the expression in question is an object of type \coq{nat} we conclude that applying any normalization algorithm to this expression we will obtain after finitely many step a numeral, i.e., we'll {\em compute} this expression. 

Consider now Martin-L{\"o}f type theory together with an added axiom \coq{A:T\_A}. While in Martin-L{\"o}f type theory strong normalization holds over any context (i.e., after the addition of any number of axioms) the canonicity theorem usually fails over most non-empty contexts. 

For example, if we obtain an object \coq{o} of type \coq{nat} using axiom \coq{funextempty} then there is no guarantee that its normal form will be a numeral or, as we say, there is guarantee that ``it will compute". However many expressions which contain an axiom will compute since the subexpressions containing the axiom get eliminated at some stage of the normalization process.

In Coq there are two main normalization algorithms which can be called by the commands \coq{Eval compute} and \coq{Eval lazy} respectively. Theoretically these algorithm are equivalent in the sense that both are supposed to always terminate and the answers produced should coincide. In practice, I have encountered many cases when \coq{Eval lazy} terminates in a reasonable amount of time while \coq{Eval compute} applied to the same expression takes too much time for me to wait it out. 

The lines in the file {\tt finitesets.v} which start with \coq{Eval compute} or \coq{Eval lazy} are tests to verify that the use of the axioms in various  proofs of finiteness does not interfere with the computability of the cardinality function \coq{fincard}. 

Note that all of the axioms which we use in this library are corollaries of the general univalence axiom. So if or when the main conjecture on constructive interpretation of the univalence will be proved we will have an algorithm which, when applied to any well-formed expression \coq{o} of type \coq{nat} which uses  any of the axioms of the library will return an expression \coq{o'} without any axioms in it and a proof that the new expression is paths-equal to \coq{o}. This algorithm will however be of a different kind than the normalization algorithms\footnote{For a recent advance in solving this problem see \url{https://github.com/simhu/cubical}.}.

\subsection{Files {\tt hz.v} and {\tt hq.v}}

In these two files we define first integers \coq{hz} and then rational numbers \coq{hq}. In both cases we follow Bourbaki approach. In the file {\tt hnat.v} we have defined a commutative rig of natural numbers. To get \coq{hz} we apply the general construction \coq{commrigtocommrng} of the ring of differences of a commutative rig from {\tt algebra1c.v}. To get \coq{hq} we apply to the integral domain \coq{hz} the general construction \coq{fldfrac} of the field of fractions from {\tt algebra1d.v}. Note that this construction requires the equality on the integral domain to be decidable. This is due to our definition \coq{isafield} of what a field is. 

At the end of both files {\tt hz.v} and {\tt hq.v} are more test computations.

\subsection{File {\tt funextfun.v}}

In this file we introduce the Univalence Axiom and prove that it implies the functional extensionality axiom \coq{funextfunax} from {\tt uu0.v} . The rest of the library does not depend on this file.

\subsection{Appendix: on the Coq system for naming and loading libraries.}

I am grateful to Dan Grayson for figuring out the answers to  many questions which I had while writing this appendix. 

At the top of the files of the Foundations library (other than \coq{uuu.v}) there are lines starting with \coq{Add LoadPath} and \coq{Require}.  These are commands which tell Coq where to look for ``libraries'' which are needed to compile the given file.  For the purpose of this explanation I will use the file \coq{uu0.v}. 

Understanding why a particular combination of these commands, the \coq{-R} options in the \coq{Makefile}, and the \coq{-R} options in the emacs variable \coq{coq-prog-args} used by the \coq{Proof General} when starting \coq{Coq} works, while a slightly different one doesn't, can be very confusing. Below I will try to describe the minimum that I believe is sufficient to understand why the particular choices made in Foundations library work as they do and to be able to predict the effect of possible small modifications of these choices.  

The \coq{Require} command in the file  \coq{uu0.v} tells Coq to load a ``library" that is called \coq{Foundations.Generalities.uuu}. The word \coq{Export}, as opposed to the word \coq{Import}, means that \coq{Foundations.Generalities.uuu} will also be loaded every time the \coq{Foundations.Generalities.uu0} (the name of the library in the file \coq{uu0.vo}) is loaded. 

Two issues contribute to the complexity of the behavior of these commands. One is how the name of the library which is contained in a given \coq{.vo} file is determined and another one is which files and directories Coq will look through when it tries to execute the \coq{Require} command and what will be the name of the library it will look for in each of these files (which will, as we will see below, be usually different from the name specified in the \coq{Require}). 

The \coq{.vo} files are created by \coq{coqc}, the non-interactive mode of Coq, using as the input a \coq{.v} file, i.e., a file which contains the humanly readable Coq code. This is what the Makefile in the top directory of Foundations library does: it calls the program \coq{coqc} for each of the \coq{.v} files in the library to produce the corresponding \coq{.vo} files. One can not, for example, experiment with \coq{uu0.v} in \coq{Proof General} until \coq{uuu.vo} has been created by running \coq{coqc} on \coq{uuu.v}. 

If no \coq{-R} option is given when \coq{coqc} command is called then the name of the library in the \coq{.vo} file created by this call is the name of the (main part of) the \coq{.v} file {\em as given to the \coq{coqc}}. In the case of \coq{uu0.v} the command \coq{coqc uu0}, run in the directory \coq{Foundations/Generalities/}, will put the name \coq{uu0} to the library in the \coq{uu0.vo} which it will produce. The command \coq{coqc Generalities/uu0} ran in the directory \coq{Foundations} will put into the \coq{uu0.vo} a ``library'' called \coq{Generalities.uu0}.

What will happen if an arbitrary \coq{-R} option is given to the \coq{coqc} command  I do not know. If the option is of the form \coq{-R "." "name"} then the name of the library in the \coq{.vo} file will be the name one would expect without the \coq{-R} option with \coq{name.} appended in front of it. The \coq{name} may itself consist of several components, e.g., it can be \coq{Foundations.Generalities}.  

Suppose now that we want to run coq on the \coq{uu0.v} file. When Coq reads the command \coq{Require Export Foundations.Generalities.uu0} it will start looking for a file whose name is \coq{uu0.vo} ``on the \coq{LoadPath}''. 

The latter expression means the following. \coq{LoadPath} is represented by a list of pairs where the second component of the pair is the actual name of a directory and the first component is an expression of the form \coq{n1.n2.\dots.nk} (where \coq{ni} are names) which Coq will use instead of the directory name internally. 

One can find the content of this list from \coq{Proof General} by  running Coq over the command \coq{Print LoadPath}.

The \coq{Add LoadPath} command adds to this list the line which you would expect from the arguments of the \coq{Add LoadPath}. A version of this command \coq{Add Rec LoadPath} will also add the lines corresponding to all of the subdirectories of the directory mentioned in the arguments (except possibly some whose names contain symbols which are not permitted in identifiers).  

If the Coq program was given \coq{-R name namedir} as an argument it will have the same effect on the \coq{LoadPath} as the command \coq{Add Rec LoadPath "namedir" name.}

When Coq encounters the line \coq{Require Export n1.n2\dots.nk.n} it does the following. First it looks for the file \coq{n.vo} in the directories \coq{dirname} which appear in \coq{LoadPath} in pair with \coq{n1.n2.\dots.nk}. It will take the first such file it finds and will check whether it contains library \coq{n1.\dots.nk.n}. If it does not it won't look for another possible match and will give an error message. If it does it will load the library. 

The content of the \coq{LoadPath} can also be modified by using \coq{-R} option when calling Coq, e.g., by customizing the variable \coq{coq-prog-args} in \coq{Proof General}. One can experiment with the results of such modifications using the \coq{Print LoadPath} command.

%\bibliography{alggeom}
%\bibliographystyle{plain}

\def\cprime{$'$}

\end{document}